\documentclass{amsart}
\usepackage{amssymb}
\newtheorem{thm}{Theorem}[subsection]
\newtheorem{prop}{Proposition}[subsection]
\newtheorem{lemma}{Lemma}[subsection]
\newtheorem{cor}{Corollary}[subsection]
\theoremstyle{remark}
\newtheorem{remark}{Remark}
\newtheorem{exam}{Example}
\newcommand{\ind}{{\rm{Index}}}
\newcommand{\la}{\langle}
\newcommand{\ra}{\rangle}

\title{Involutive equivalence bimodules and inclusions of $C^{*}$-algebras
with Watatani index 2}
\author{Kazunori Kodaka \and Tamotsu Teruya}
\address{Department of Mathematical Sciences, Faculty of Science,
Ryukyu University \\
Nishihara-cho, Okinawa 903-0213, Japan}
\address{\sl{E-mail address}: \rm{kodaka@math.u-ryukyu.ac.jp}}
\address{\sl{E-mail address}: \rm{teruya@math.u-ryukyu.ac.jp}}

\begin{document}

\maketitle
\baselineskip=20pt
\centerline{Abstract}
\bigskip
Let $A$ be a unital $C^*$-algebra.
We shall introduce involutive $A$-$A$-equivalence bimodules and prove that
any $C^*$-algebra containing $A$ with Watatani index 2 is constructed
by an involutive $A$-$A$-equivalence bimodule.

\bigskip

AMS 2000 Mathematical Subject Classification: Primary 46L08, Secondary
46L40.


\section{introduction}\label{intro}
V. Jones introduced an index theory for $II_1$ factors in \cite{Jones:factor}.
One of his motivations is the Goldman's theorem, which says that if $M$ is a type
$II_1 $ factor and $N\subset M$ is a subfactor with the Jones index
$[M:N]=2$, then there is a crossed
product decomposition $M = N \rtimes_{\alpha}{\Bbb Z}_2$, where $
{\Bbb Z}_2$ is the group ${\Bbb Z}/{2\Bbb Z}$ of order two.		
Since Jones index theory is extended to $C^*$-algebras by Y. Watatani,
it is worth to investigate Goldman type theorem for inclusions of simple
$C^*$-algebras. 
In the present paper, we shall study the inclusion $A \subset B$ of
$C^*$-algebras
with a conditional expectation $E:B \to A$ of Index $E =2$.
In 4.2 Examples, we shall show that Goldman type theorem does not hold for
inclusions of simple $C^*$-algebras in general by exhibiting examples
of inclusions like a non-commutative sphere in an irrational rotation $C^*$-algebra
$A_{\theta}$ and irrational rotation $C^*$-algebras
$A_{2\theta} \subset A_{\theta}$
with different angles. Therefore there occurs the following natural question:
What kind of structures are there in the inclusion of $C^*$-algebras with index 2? 
We shall answer the question in the present paper:
Any inclusion of $C^*$-algebras with index two
gives an involutive equivalence bimodule.
	
Let us explain the notion of involutive 
equivalence bimodules. Consider a typical situation, that is, the inclusion $A \subset B$ is given by the crossed product
$B= A \rtimes_{\alpha}{\Bbb Z}_2$ by some action $\alpha : {\Bbb Z}_2 \to Aut (A)$.
Then the canonical conditional expectation $E:B \to A$ has Index $E =2$. 
Moreover there exists the dual action $\hat{\alpha}:{\Bbb Z}_2 \to Aut (B)$ such that
$$
(A \rtimes_{\alpha}{\Bbb Z}_2)\rtimes_{\hat{\alpha}}{\Bbb Z}_2 
\cong A \otimes M_2(\Bbb C) ,
$$
where $M_2 (\Bbb C )$ is the $2\times 2$-matrix algebra over $\Bbb C$.
It is well known that the $C^*$-basic construction $C^*\langle B, e_A \rangle$ is exactly 
$(A \rtimes_{\alpha}{\Bbb Z}_2)\rtimes_{\hat{\alpha}}{\Bbb Z}_2$. 
Then the Jones projection $e_A$ corresponds to
the projection $e_{11} = diag(1,0)$ 
and $ 1 - e_A$ corresponds to $e_{22}= diag(0, 1)$, where
$diag(\lambda, \mu)$ is a $2\times2$-diagonal matrix with
diagonal elements $\lambda$, $\mu$. Let
$X= e_{11}(A \otimes M_2(\Bbb C))e_{22}$. Then $X$ is an
$A$-$A$-equivalence bimodule in the natural way.
There exists a natural involution on $X$ such that 
$$
x^\sharp = 
\left(\begin{matrix}
0 & z^* \\
0 & 0
\end{matrix}\right)
\quad \text{for} \quad
x =
\left(\begin{matrix}
0 & z \\
0 & 0
\end{matrix}\right).
$$	
We pick up these properties to define the notion of 
involutive equivalence bimodules.
In Theorem 3.3.1, we shall show that even if $B$ is not a crossed product of $A$, 
the inclusion of $C^*$-algebras with index 2 gives an involutive 
$A$-$A$-equivalence bimodule. Moreover the set of inclusions of  
$C^*$-algebras with index 2 has a one to one correspondence with the set 
of involutive $A$-$A$-equivalence bimodules up to isomorphisms.
	
In Proposition 4.1.2, we shall characterize the subclass such that $B$ is the 
twisted crossed product of $A$ by a partially inner $C^*$-dynamical system 
studied by Green, Olsen and Pedersen. The characterization is given by 
the von Neumann equivalence of $e_A$ and $1-e_A$ in $C^*\langle B, e_A \rangle$. 

\section{Preliminaries}\label{sec:pre}

\subsection{Some results for inclusions with index 2}

Let $B$ be a unital $C^*$-algebra and $A$ a $C^*$-subalgebra of $B$ with a
common unit.
Let $E$ be a conditional expectation of $B$ onto $A$ with $ 1 <\ind E<
\infty$.
Then by Watatani \cite{Watatani:index} we have the $C^*$-basic construction
$C^*\la B, e_A \ra$ where $e_A$ is the Jones projection induced by $E$.
Let $\widetilde{E}$ be the dual conditional expectation of $C^*\la B, e_A
\ra$
onto $B$ defined by
$$
\widetilde{E}(ae_A b) = \dfrac{1}{t}ab
\quad
\text{for any}\quad  a, b \in B,
$$ where $t = \ind E$.
Let $F$ be a linear map of
$(1-e_A)C^*\la B, e_A \ra(1-e_A)$ to $A(1-e_A)$ defined by
$$
F(a)=\dfrac{t}{t-1}(E \circ \widetilde{E})(a)(1-e_A)
$$
for any $ a \in (1-e_A)C^*\la B, e_A \ra(1-e_A)$.
By routine computations we can see that $F$ is a
conditional expectation of
$(1-e_A)C^*\la B, e_A \ra(1-e_A)$ onto $A(1-e_A)$.

\begin{lemma}\label{lemma:local}
With the above notations, let $\{(x_i, x_i^*)\}_{i=1}^n$ be a quasi-basis
for $E$. Then
$$
\{\sqrt{t-1}(1-e_A)x_je_Ax_i(1-e_A), \sqrt{t-1}(1-e_A)x_i^*e_Ax_j^*(1-e_A)
\}_{i,j=1}^n
$$
is a quasi-basis for $F$. Furthermore $\ind F =(t-1)^2(1-e_A)$.
\end{lemma}
\begin{proof}
This is immediate by direct computations.
\end{proof}

\begin{cor}\label{cor:cut}
We suppose that
$\ind E =2$. Then
$$
(1-e_A)C^*\la B, e_A \ra (1-e_A) = A(1-e_A) \cong A.
$$
\end{cor}
\begin{proof}
By Lemma \ref{lemma:local} there is a conditional expectation $F$ of
$(1-e_A)C^*\la B, e_A \ra (1-e_A)$ onto $A(1-e_A)$ and
$$
\ind F = (\ind E - 1)^2(1-e_A).
$$
Since $\ind E =2$, $\ind F = 1-e_A$. Hence by Watatani
\cite{Watatani:index},
$$
(1-e_A)C^*\la B, e_A \ra (1-e_A)=A(1-e_A).
$$
If $a(1-e_A) = 0$, for $a\in A$, then $a=2\widetilde{E}(a(1-e_A))=0$.
Therefore the map $a \mapsto a(1-e_A)$ is injective.
And hence $A(1-e_A) \cong A$ as desired.
\end{proof}

\begin{lemma}\label{lemma:cut}
With the same assumptions as in Lemma \ref{lemma:local}, we suppose that
$\ind E =2$. Then
for any $b \in B$,
$$
(1-e_A)b(1-e_A) = E(b)(1-e_A).
$$
\end{lemma}
\begin{proof}
By Corollary \ref{cor:cut} there exists $a \in A$ such that
$(1-e_A)b(1-e_A)= a(1-e_A)$.
Therefore $a= 2\widetilde{E}(a(1-e_A))
= 2\widetilde{E}((1-e_A)b(1-e_A)) = E(b)$.
This completes the proof.
\end{proof}

\begin{prop}\label{prop:aut}
With the same assumptions as in Lemma \ref{lemma:local}, we suppose that
$\ind E =2$.  Then  there is a unitary element
$U \in C^*\la B, e_A\ra$ satisfying the following conditions:
\begin{enumerate}
\item $U^2 =1$,
\item $UbU^*=2E(b)-b$ for $b \in B$.
\end{enumerate}
Hence if we denote by $\beta$ the restriction of $Ad(U)$ to
$B$, $\beta$ is an automorphism of $B$ with
$\beta^2= id$ and $B^{\beta}=A$.
\end{prop}

\begin{proof}
By Lemma \ref{lemma:cut}, for any $b \in B$
$$
(1-e_A)b(1-e_A)=E(b)(1-e_A) = E(b)-E(b)e_A .
$$
On the other hand
$$
(1-e_A)b(1-e_A)
= b -e_Ab -be_A + E(b)e_A .
$$
Therefore
$$
E(b) =  b -e_Ab -be_A + 2E(b)e_A.
$$
Let $U$ be a unitary element defined by $U = 2e_A -1$.
Then by the above equation for any $b \in B$
\begin{align*}
UbU^*
&= 2(b -e_Ab -be_A + 2E(b)e_A) - b \\
&= 2E(b)-b.
\end{align*}
Thus we obtain the conclusion.
\end{proof}

\begin{remark}\label{rem:cross}
By the above proposition, $E(b)=\frac{1}{2}(b+\beta(b))$.
\end{remark}

\begin{lemma}\label{lem:equi}
Let $B$ be a unital $C^*$-algebra and $A$ a $C^*$-subalgebra of $B$ with a
common unit.
Let $E$ be a conditional expectation of $B$ onto $A$ with $\ind E=2$.
Then we have
$$C^*\la B,e_A \ra \cong B \rtimes_{\beta} \mathbb{Z}_2.$$
\end{lemma}

\begin{proof}
We may assume that
$B \rtimes_{\beta}\mathbb{Z}_2$
acts on the Hilbert space
$l^2(\mathbb{Z}_2,H)$ faithfully, where $H$ is some Hilbert space on
which $B$ acts faithfully. Let $W$ be a unitary element in
$B \rtimes_{\beta}\mathbb{Z}_2$ with $\beta = Ad(W)$, $W^2=1$.
Let $e=\frac{1}{2}(W+1)$. Then $e$ is a projection in
$B \rtimes_{\beta}\mathbb{Z}_2$ and
$ebe=E(b)e$ for any $b \in B$.
In fact,
$$
ebe= \frac{1}{4}(WbW+bW+Wb+b).
$$
On the other hand by Remark \ref{rem:cross},
$$
E(b)e = \frac{1}{2}(b+\beta(b))\frac{1}{2}(W+1) =
\frac{1}{4}(WbW+bW+Wb+b).
$$
Hence $ebe=E(b)e$
for
$b\in B$.
Also
$A \ni a \mapsto ae \in B \rtimes_{\beta}\mathbb{Z}_2$ is
injective.
In fact, if
$ae=0$,
$aW+a=0$. Let
$ \widehat{\beta}$
be the dual action of
$\beta$.
Then $0= \widehat{\beta}(aW + a) = -aW+ a$. Thus $2a =0$, i.e., $a=0$.
Hence by Watatani \cite [Proposition 2.2.11]{Watatani:index},
$C^*\la B, e_A \ra \cong B \rtimes_{\beta}\mathbb{Z}_2$.
\end{proof}
\begin{remark}\label{rem:kappa}

\begin{enumerate}

        \item  By the proofs of Watatani \cite[Propositions 2.2.7 and
2.2.11]{Watatani:index},
             we see
                that $\kappa(b) =b$ for any $b \in B$, where $\kappa$
is the isomorphism
of
                $C^*\la B,e_A \ra$ onto $B \rtimes_{\beta} \mathbb{Z}_2$
in Lemma
\ref{lem:equi}.
        \item The above lemma is obtained in Kajiwara and Watatani
                \cite[Theorem 5.13]{KW}
\end{enumerate}
\end{remark}

By Lemma \ref{lem:equi} and Remark \ref{rem:kappa}, we regard
$\widehat{\beta}$ as an automorphism of $C^*\la B,e_A \ra $ with
$\widehat{\beta}(b)=b$ for any $b\in B, {\widehat{\beta}}^2=id$ and
$\widehat{\beta}(e_A) = 1 - e_A$.

\begin{lemma}\label{lem:fix}
With the same assumptions as in Lemma \ref{lem:equi},
$$
{C^*\la B,e_A \ra}^{\widehat{\beta}} = B.
$$
\end{lemma}
\begin{proof}
By Lemma \ref{lem:equi} for any $x \in C^*\la B,e_A \ra$, we can write
$x = b_1 + b_2U$, where $b_1, b_2 \in B$. We suppose that
$\widehat{\beta}(x) =x$.
Then $b_1 - b_2U = b_1 + b_2U$.
Thus $b_2 =0$.
Hence $x= b_1 \in B$.
Since it is clear that $B \subset {C^*\la B,e_A \ra}^{\widehat{\beta}}$,
the lemma is proved
\end{proof}

\subsection{Involutive equivalence bimodules}

Let $A$ be a unital $C^*$-algebra and $X (={}_AX_A)$ an  $A$-$A$-
equivalence
bimodule.
$X$ is {\it involutive} if there exists a conjugate linear map $x \mapsto
x^\sharp$ on $X$,
such that
\begin{enumerate}
\item $(x^\sharp)^\sharp = x$, \quad $ x \in X$,
\item $(a\cdot x \cdot b)^\sharp =
b^* \cdot x^\sharp \cdot a^*$, \quad $x\in X$, $a,b
\in A$,
\item ${}_A\la x, y^\sharp \ra = \la x^\sharp, y\ra_A$, \quad $x,y \in X$,
\end{enumerate}
where ${}_A\la,\ra$ and $\la, \ra_A$ are the left and the
right $A$-valued inner
products on
$X$, respectively.
We call the above conjugate linear map an involution on $X$.

For an $A$-$A$-equvalence bimodule $X$, we define its
dual bimodule. Let $\widetilde{X}$ be $X$ itself when it is considered
as a set. We write $\widetilde{x}$ when $x$ is considered in
$\widetilde{X}$. $\widetilde{X}$ is made into an equivalence
$A$-$A$-bimodule as follows:
\begin{enumerate}
\item $\widetilde{x}+\widetilde{y}=\widetilde{x+y}$,
$\lambda\widetilde{x}=\widetilde{\bar{\lambda}x}$
for any $x, y\in X$ and $\lambda\in\Bbb{C}$,
\item $b \cdot\widetilde{x}\cdot a=\widetilde{a^* \cdot x\cdot b^*}$
for any $a, b\in A$ and $x\in X$,
\item ${}_A \la \widetilde{x}, \widetilde{y} \ra=\la x, y\ra_A$,
$\la\widetilde{x}, \widetilde{y} \ra_A ={}_A \la x, y \ra$
for any $x, y\in X$.
\end{enumerate}

\begin{lemma}\label{V:iso}
Let $V$ be a map of an involutive $A$-$A$-equivalence
bimodule $X$ onto its dual bimodule $\widetilde{X}$ defined by
$V(x)= \widetilde{x^\sharp}$, where $\widetilde{x}$ means $x$ as viewed
as an element in $\widetilde{X}$. Then $V$ is
an $A$-$A$-equivalence bimodule isomorphism of $X$ onto $\widetilde{X}$.
\end{lemma}
\begin{proof}
This is immediate by routine computations.
\end{proof}

%
%

\section{Correspondence between involutive equivalence bimodules and
inclusions of
$C^*$-algebras with index 2}\label{subsection:const1}

Let $A$ be a unital $C^*$-algebra and we denote by $(B,E)$ a pair of
a unital $C^*$-algebra $B$ including $A$ as a $C^*$-subalgebra of $B$
with a common unit and a
conditional
expectation $E$ of $B$ onto $A$ with \ind $E=2$.
Let $\mathcal{L}$ be the set of all such pairs $(B, E)$ as above.
We define an equivalence relation $\sim$ in $\mathcal{L}$ as follows:
for $(B,E), (B_1,E_1) \in \mathcal{L}$, $(B,E)\sim(B_1,E_1)$
if and only if there is an isomorphism $\pi$ of $B$ onto $B_1$
such that $\pi(a) = a$ for any $a \in A$ and
$E_1 \circ \pi = E$. We denote by $[B,E]$ the equivalence class of $(B,E)$.

Let $\mathcal{M}$ be the set of all involutive $A$-$A$-equivalence
bimodules.
We define an equivalence relation $\sim$ in $\mathcal{M}$ as follows:
for $X, Y \in \mathcal{M}$,
$X\sim Y$ if and only if there is an $A$-$A$-equivalence
bimodule isomorphism
$\rho$ of $X$ onto $Y$ with $\rho(x^\sharp) = \rho(x)^\sharp$.
We call $\rho$ an involutive $A$-$A$-equivalence
bimodule isomorphism of $X$ onto $Y$.
We denote by $[X]$ the equivalence class of $X$.

\subsection{Construction of a map from $\mathcal{L}/\sim$ to
$\mathcal{M}/\sim$}

We shall use the same notations as in section \ref{sec:pre}.

Let $B$ be a unital $C^*$-algebra and $A$ a $C^*$-subalgebra of $B$ with a
common unit.
Let $E$ be a conditional expectation of $B$ onto $A$ with $\ind E=2$.
Then by Watatani \cite{Watatani:index} and Corollary \ref{cor:cut}
\begin{enumerate}
\item $e_AC^*\la B,e_A \ra e_A = Ae_A \cong A$,
\item $(1-e_A)C^*\la B,e_A \ra (1-e_A) = A(1-e_A) \cong A$.
\end{enumerate}
Let $\psi$ be an isomorphism of $A$ onto $Ae_A$ defined by $\psi(a)=ae_A$
for any $a \in A$ and $\phi$ an isomorphism of $A$ onto $A(1-e_A)$
defined by $\phi = \widehat{\beta}\circ \psi$.
Let $X_{(B,E)}=X_B =e_AC^*\la B,e_A \ra (1-e_A)$. We regard $X_B$ as a
Hilbert $A$-$A$-bimodule in the following way:
for any $a,b \in A$ and $x \in X_B$,
$a\cdot x \cdot b = \psi(a)x\phi(b) = axb$.
For any $x,y \in X_B$,
${}_{A}\la x , y\ra = \psi^{-1}(xy^*)$,
$\la x,  y\ra_A = \phi^{-1}(x^*y)$.

\begin{lemma}\label{lem:bimo}
With the above notations, $X_B$ is an $A$-$A$-equivalence bimodule.
\end{lemma}
\begin{proof}
This is immediate by routine computations.
\end{proof}

Let $x \mapsto x^\sharp$ be a conjugate linear map on $X_B$
defined by $x^\sharp = \widehat{\beta}(x^*)$
for any $x \in X_B$.
Since $\widehat{\beta}^2 = id$, $(x^\sharp)^\sharp = x$.
Since $\widehat{\beta}(a) =a $ for any $a \in A$,
$(a\cdot x \cdot b)^\sharp = \widehat{\beta}( b^* x^* a^*)
= b^* \cdot x^\sharp \cdot a^*$ for $x\in X$, $a,b \in A$.
Furthermore, for $x, y \in X_B$
${}_A\la x, y^\sharp \ra=\la x^{\sharp}, y \ra_A $
by an easy calculation.
Therefore $X_B$ is an element in $\mathcal{M}$.

\begin{remark}
$\widetilde{X_B}$ is isomorphic to $(1-e_A)C^*\la B,e_A \ra e_A$
as $A$-$A$-equivalence bimodules.
Indeed, the map
$(1-e_A)C^*\la B,e_A \ra e_A \ni (1-e_A)xe_A \mapsto
\widetilde{e_Ax^*(1-e_A)}$,
$x \in C^*\la B,e_A \ra$ gives an $A$-$A$-equivalence
bimodule isomorphism of
$(1-e_A)C^*\la B,e_A \ra e_A$ onto $\widetilde{X_B}$,
where $\widetilde{y}$ means $y$ viewed as an element in $\widetilde{X_B}$
for any $y \in X_B$.
Sometimes, we identify $\widetilde{X_B}$ with $(1-e_A)C^*\la B,e_A \ra e_A$.
\end{remark}

Let $\mathcal{F}$ be a map from $\mathcal{L}/\sim$ to $\mathcal{M}/\sim$
defined by $\mathcal{F}([B,E]) = [X_B]$ for any $[B,E] \in
\mathcal{L}/\sim$.

\begin{lemma}
With the above notations, $\mathcal{F}$ is well-defined.
\end{lemma}
\begin{proof}
Let $(B,E), (B_1,E_1) \in \mathcal{L}$ with
$(B,E)\sim (B_1,E_1)$.
Let
$X_B$ and
$X_{B_1}$ be elements in
$\mathcal{M}$
defined by $(B,E)$ and $(B_1,E_1)$, respectively.
Since $(B,E)\sim(B_1,E_1)$, there is an isomorphism $\pi$ of $B$
onto $B_1$ such that $\pi(a)= a$ for any $a \in A$ and $E_1\circ \pi =E$.
Let $\widetilde{\pi}$ be a homomorphism of the linear span of
$\{\, be_Ac, |\, b, c \in B\,\} $ to $C^*\la B_1, e_{A,1}\ra$ defined by
$\widetilde{\pi}(be_Ac) = \pi(b)e_{A,1}\pi(c)$ for any
$b,c \in B$.
Then for $b_i, c_i \in B (i=1,2, \cdots, n)$ and $a \in B$,
\begin{align*}
&\left\| \widetilde{\pi}\left( \sum_{i=1}^{n}b_ie_Ac_i \right)
\pi(a)\right\|^2
=\left\| \sum_{i=1}^{n}\pi(b_i)E_1(\pi(c_ia))\right\|^2 \\
&=\left\|
\sum_{i,j=1}^{n}E_1(\pi(a^*c_i^*))E_1(\pi(b_i^*b_j))E_1(\pi(c_ja))\right\| \\
& =\left\| \sum_{i,j=1}^{n}E(a^*c_i^*)E(b_i^*b_j)E(c_ja)\right\|.
\end{align*}
On the other hand
$$
\left\| \sum_{i=1}^{n}b_ie_Ac_ia\right\|^2
=\left\| \sum_{i=1}^{n}b_iE(c_ia)\right\|^2
=\left\| \sum_{i,j=1}^{n}E(a^*c_i^*)E(b_i^*b_j)E(c_ja)\right\|.
$$
Hence
\begin{align*}
\left\| \widetilde{\pi}\left( \sum_{i=1}^{n}b_ie_Ac_i \right)\right\|
&= \sup \left\{\left\|\widetilde{\pi}\left( \sum_{i=1}^nb_ie_Ac_i
\right)\pi(a)\right\| |
\left\|E_1(\pi(a)^*\pi(a))\right\| =1, a\in B \right\} \\
&= \sup \left\{ \left\| \sum_{i=1}^nb_ie_Ac_i a \right\| |
\left\|E(a^*a) \right\| =1, a\in B \right\} \\
&=  \left\| \sum_{i=1}^nb_ie_Ac_i \right\|.
\end{align*}
Thus $\widetilde{\pi}$ can be extended to an isomorphism of
$C^*\la B, e_A\ra$ onto $C^*\la B_1, e_{A,1}\ra$.
Hence $\widetilde{\pi}$ is an involutive $A$-$A$-equivalence
bimodule isomorphism of
$X_B$ onto $X_{B_1}$
since
$\widetilde{\pi}(e_A)=e_{A,1}$.
In fact, for $a \in A$ and $x \in C^*\la B, e_A\ra $
$$
\widetilde{\pi}(a\cdot e_Ax(1-e_A))
= e_{A,1}a\pi(x)(1-e_{A,1})
= a \cdot \widetilde{\pi}(e_Ax(1-e_A)),
$$
Similarly
$$
\widetilde{\pi}(e_Ax(1-e_A)\cdot a)
= \widetilde{\pi}(e_Ax(1-e_A))\cdot a.
$$
Also, for $x,y \in C^*\la B, e_A\ra$,
\begin{align*}
{}_A\la \widetilde{\pi}(e_Ax(1-e_A)), \widetilde{\pi}(e_Ay(1-e_A))\ra
&= \left(\psi_1^{-1}\circ \widetilde{\pi} \right)(e_Ax(1-e_A)y^*e_A) \\
&= {}_A\la e_Ax(1-e_A), e_Ay(1-e_A)\ra,\\
& {} \\
\la \widetilde{\pi}(e_Ax(1-e_A)), \widetilde{\pi}(e_Ay(1-e_A))\ra_A
&= \phi^{-1}((1-e_A)x^*e_Ay(1-e_A)) \\
&= \la e_Ax(1-e_A), e_Ay(1-e_A)\ra_A
\end{align*}
since $\psi_1^{-1} = \widetilde{\pi} \circ \psi$ and
$\widetilde{\pi}\circ\widehat{\beta} = \widehat{\beta}_1\circ
\widetilde{\pi}$.
Furthermore, for any $x \in C^* \la B,e_A\ra$
\begin{align*}
\widetilde{\pi}((e_Ax(1-e_A))^\sharp)
&= \widetilde{\pi}(e_A\widehat{\beta}(x)^*(1-e_A)) \\
&= (e_{A,1}\widetilde{\pi}(x)(1-e_{A,1}))^\sharp
= \widetilde{\pi}(e_Ax(1-e_A))^\sharp.
\end{align*}
Therefore $X_B \sim X_{B_1}$ in $ \mathcal{M}$.
\end{proof}

%
%

\subsection{Construction of a map from $\mathcal{M}/\sim$ to
$\mathcal{L}/\sim$}
Let $X \in \mathcal{M}$. Following Brown, Green and Rieffel
\cite{BGR:stable},
we can define the linking algebra $L$ for an $A$-$A$-equivalence
bimodule $X$. Let
$$
L_0= \left\{
\begin{bmatrix}
 a & x \\
\widetilde{y} & b
\end{bmatrix} \quad
| \quad a, b \in A, x,y \in X \right\},
$$
where $\widetilde{y}$ means $y$ viewed as 
an element in the dual bimodule
$\widetilde{X}$ of $X$.
In the same way as in Brown, Green and Rieffel \cite{BGR:stable}
we can see that $L_0$ is a $*$-algebra.
Also we regard $L_0$ as a $*$-subalgebra acting on the
right Hilbert $A$-module $X \oplus A$.
Hence we can define an operator norm in $L_0$ acting on $X \oplus A$.
We define $L$ as the above operator norm closure of $L_0$.
But, since $X$ is complete, in this case $ L= L_0^- = L_0$.
Let $B_X$ be a subset of $L$ defined by
$$
B_X =
\left\{
\begin{bmatrix}
a & x \\
\widetilde{x^\sharp} & a
\end{bmatrix}
\quad
|
\quad
 a \in A, x \in X
\right\}.
$$
By direct computations, we can see that
$B_X$ is a $*$-subalgebra of $L$ and since $X$ is complete, $B_X$ is
closed in $L$, that is, $B_X$ is a $C^*$-subalgebra of $L$.
We regard $A$ as a $C^*$-subalgebra
$
\left\{
\begin{bmatrix}
a & 0 \\
0 & a
\end{bmatrix} \ |
\ a \in A
\right\}
$
of $B_X$.
Let $E_X$ be a linear map of $B_X$ onto $A$ defined by
$
E_X\left(
\begin{bmatrix}
a & x \\
\widetilde{x^\sharp} & a
\end{bmatrix}
\right)
=
\begin{bmatrix}
a&0\\
0&a
\end{bmatrix}
$
for any
$
\begin{bmatrix}
a & x \\
\widetilde{x^\sharp} & a
\end{bmatrix}
\in B_X
$.
Then by easy computations $E_X$ is a conditional expectation
of $B_X $ onto $A$.

\begin{lemma}
With the above notations, $\ind E_X = 2$.
\end{lemma}
\begin{proof}
There are elements $z_1, \dots, z_n, \, y_1, \dots, y_n \in X$
such that
$\sum_{i=1}^n \la z_i, y_i\ra_A= 1$ by
Rieffel \cite[the proof of Proposition 2.1]{Reiffel:irr}
since $X$ is an $A$-$A$-equivalence bimodule.
For $i = 1,2,\cdots,n$ let $w_i$ be an element in $X$
with $w_i = z_i^\sharp$.
Then
$
\left\{
\left(
\begin{bmatrix}
1&0\\
0&1
\end{bmatrix},
\begin{bmatrix}
1&0\\
0&1
\end{bmatrix}
\right)
\right\}
\cup
\left\{
\left(
\begin{bmatrix}
0&w_i \\
\widetilde{w_i^\sharp} & 0
\end{bmatrix},
\begin{bmatrix}
0&y_i \\
\widetilde{y_i^\sharp} & 0
\end{bmatrix}
\right)
| \, i=1,2, \cdots, n
\right\}
$
is a quasi-basis for $E_X$ by direct computations.
In fact, for
$
\begin{bmatrix}
a & x \\
\widetilde{x^\sharp} & a
\end{bmatrix} \in B_X
$
\begin{align*}
&E\left(
\begin{bmatrix}
a&x \\ \widetilde{x^\sharp}&a
\end{bmatrix}
\begin{bmatrix}
1&0 \\ 0&1
\end{bmatrix}
\right)
\begin{bmatrix}
1&0 \\ 0&1
\end{bmatrix}
=
\begin{bmatrix}
a&0 \\ 0&a
\end{bmatrix}, \\
&E\left(
\begin{bmatrix}
a&x \\ \widetilde{x^\sharp}&a
\end{bmatrix}
\begin{bmatrix}
0&w_i \\ \widetilde{w_i^\sharp}&0
\end{bmatrix}
\right)
\begin{bmatrix}
0&y_i \\ \widetilde{y_i^\sharp}&0
\end{bmatrix}
=
\begin{bmatrix}
0 & {}_A\la x, w_i^\sharp\ra y_i \\
\la x^\sharp , w_i\ra_A \widetilde{y_i^\sharp} & 0
\end{bmatrix}.
\end{align*}
Also,
$$
\sum_{i=1}^n {}_A\la x, w_i^\sharp \ra y_i
= \sum_{i=1}^n x \la w_i^\sharp , y_i \ra_A
= x,
$$
$$
\sum_{i=1}^n \la x^\sharp, w_i \ra_A \widetilde{y_i^\sharp}
= \sum_{i=1}^n {}_A\la x, w_i^\sharp \ra \widetilde{y_i^\sharp}
= \sum_{i=1}^n V({}_A\la x, w_i^\sharp \ra y_i)
= \widetilde{x^\sharp},
$$
where $V$ is an $A$-$A$-equivalence bimodule isomorphism
defined in Lemma \ref{V:iso}.
Hence
$$
E\left(
\begin{bmatrix}
a & x \\ \widetilde{x^\sharp}& a
\end{bmatrix}
\begin{bmatrix}
1 & 0 \\ 0 & 1
\end{bmatrix}
\right)
\begin{bmatrix}
1 & 0 \\ 0 & 1
\end{bmatrix} + \sum_{i=1}^n E\left(
\begin{bmatrix}
a & x \\ \widetilde{x^\sharp}& a
\end{bmatrix}
\begin{bmatrix}
0 & w_i \\ \widetilde{w_i^\sharp} & 0
\end{bmatrix}
\right)
\begin{bmatrix}
0 & y_i \\
\widetilde{y_i^\sharp} & 0
\end{bmatrix}
= \begin{bmatrix}
a & x \\ \widetilde{x^\sharp}& a
\end{bmatrix}.
$$
Similarly
$$
\begin{bmatrix}
1 & 0 \\ 0 & 1
\end{bmatrix}
E\left(
\begin{bmatrix}
1 & 0 \\ 0 & 1
\end{bmatrix}
\begin{bmatrix}
a & x \\ \widetilde{x^\sharp}& a
\end{bmatrix}
\right) +
\sum_{i=1}^n
\begin{bmatrix}
0& w_i \\ \widetilde{w_i^\sharp} & 0
\end{bmatrix}
E \left(
\begin{bmatrix}
0 & y_i \\ \widetilde{y_i^\sharp} & 0
\end{bmatrix}
\begin{bmatrix}
a & x \\ \widetilde{x^\sharp} & a
\end{bmatrix}
\right) =
\begin{bmatrix}
a & x \\ \widetilde{x^\sharp} & a
\end{bmatrix}.
$$
Thus
$$
\ind E_X
=
\begin{bmatrix} 1 & 0 \\ 0 &1 \end{bmatrix}
+ \sum_{i=1}^n \begin{bmatrix} 0 & w_i \\ \widetilde{w_i^\sharp}& 0
\end{bmatrix}
\begin{bmatrix} 0 & y_i \\ \widetilde{y_i^\sharp} & 0 \end{bmatrix}
= \begin{bmatrix} 2 & 0 \\ 0 & 2 \end{bmatrix}.
$$
Therefore we obtain the conclusion.
\end{proof}

\begin{remark}\label{remark:basic constraction}
Let $e$ be an element in $L(=L_0)$ defined by
$\begin{bmatrix}
1 & 0 \\
0 & 0
\end{bmatrix}$.
Then it is obvious that for any $ b \in B_X, \ ebe = E_X(b)e$.
Furthermore the map
$\begin{bmatrix}
a & 0 \\
0 & a
\end{bmatrix}
\mapsto
e \begin{bmatrix}
a & 0 \\
0 & a
\end{bmatrix}
= \begin{bmatrix}
a & 0 \\
0 & 0
\end{bmatrix}$
for $ a \in A $ is injective.
And hence $ L $ is the $C^*$-basic construction of $ A \subset B $ by
Watatani \cite{Watatani:index}.

\end{remark}

Let $\mathcal{G}$ be a map from $\mathcal{M}/\sim$ to $\mathcal{L}/\sim$
defined by
$\mathcal{G}([X]) = [B_X, E_X]$ for any $[X] \in \mathcal{M}/\sim$.

\begin{lemma} $\mathcal{G}$ is well-defined.
\end{lemma}
\begin{proof}
Let $X, X_1 \in \mathcal{M}$ with $X\sim X_1$.
Let $(B_X,E_X)$ and $(B_{X_1},E_{X_1})$ be elements in $\mathcal{L}$
induced by  $X$ and $X_1$, respectively.
Since $X \sim X_1$, there is an involutive $A$-$A$-equivalence
bimodule isomorphism
$\rho$ of $X$ onto $X_1$. Let $\pi$ be a map of $B_X$ to $B_{X_1}$ defined
by for
any
$\begin{bmatrix} a&x \\ \widetilde{x^\sharp} & a \end{bmatrix} \in B_X$,
$\pi \left(
\begin{bmatrix} a&x \\ \widetilde{x^\sharp} & a \end{bmatrix}\right)
= \begin{bmatrix} a&\rho(x) \\ \widetilde{\rho(x)^\sharp} & a \end{bmatrix}$.
Then it is clear that $\pi$ is linear.
For $\begin{bmatrix} a&x \\ \widetilde{x^\sharp}&a \end{bmatrix} \in B_X$,
$$
\pi \left( \begin{bmatrix} a&x \\ \widetilde{x^\sharp} & a
\end{bmatrix}\right)^*
= \begin{bmatrix} a&\rho(x) \\ \widetilde{\rho(x)^\sharp} & a
\end{bmatrix}^*
= \begin{bmatrix} a^*& \rho(x^\sharp)
\\ \widetilde{\rho(x)} & a^* \end{bmatrix}
= \pi \left( \begin{bmatrix} a&x \\ \widetilde{x^\sharp} & a
\end{bmatrix}^* \right).
$$
Also for $\begin{bmatrix} a & x \\ \widetilde{x^\sharp} & a \end{bmatrix}$
and $\begin{bmatrix}b & y \\ \widetilde{y^\sharp} & b \end{bmatrix} \in
B_X$,
$$
\pi \left( \begin{bmatrix} a & x\\ \widetilde{x^\sharp} & a \end{bmatrix}
\begin{bmatrix}b & y \\ \widetilde{y^\sharp} & b \end{bmatrix}\right)
=
\begin{bmatrix} ab+ {}_A \la x, y^\sharp \ra & \rho(ay + xb) \\
\widetilde{\rho(xb + ay)^\sharp} & \la x^\sharp, y \ra_A + ab \end{bmatrix},
$$
and
\begin{align*}
\pi \left( \begin{bmatrix} a & x\\ \widetilde{x^\sharp} & a \end{bmatrix}
\right)
&\pi \left(\begin{bmatrix}b & y \\ \widetilde{y^\sharp} & b
\end{bmatrix}\right)
= \begin{bmatrix} ab + {}_A \la \rho(x), \rho(y^\sharp) \ra &
\rho(ay + xb) \\ \widetilde{\rho(xb+ ay)^\sharp} &
\la \rho(x^\sharp), \rho(y) \ra_A + ab \end{bmatrix}\\
&= \begin{bmatrix} ab+ {}_A \la x, y^\sharp \ra & \rho(ay + xb) \\
\widetilde{\rho(xb + ay)^\sharp} & \la x^\sharp, y \ra_A + ab \end{bmatrix}
=\pi \left( \begin{bmatrix} a & x\\ \widetilde{x^\sharp} & a \end{bmatrix}
\begin{bmatrix}b & y \\ \widetilde{y^\sharp} & b \end{bmatrix}\right).
\end{align*}
Hence $\pi$ is a homomorphism of $B_X$ to $B_{X_1}$.
Furthermore, by the definition of $\pi$, $\pi$ is a bijection and
$\pi \left(\begin{bmatrix}a & 0 \\ 0 & a \end{bmatrix}\right) =
\begin{bmatrix}a & 0 \\ 0 & a \end{bmatrix}$ for any $a \in A$.
And for $\begin{bmatrix} a & x \\ \widetilde{x^\sharp} & a \end{bmatrix} \in
B_X$
$$
(E_1 \circ \pi)\left( \begin{bmatrix} a & x \\ \widetilde{x^\sharp} & a
\end{bmatrix}\right)
= E_1 \left(  \begin{bmatrix} a & \rho(x) \\ \widetilde{\rho(x)^\sharp} & a
\end{bmatrix}\right)
= \begin{bmatrix} a & 0 \\ 0&a \end{bmatrix}
= E \left( \begin{bmatrix} a & x \\ \widetilde{x^\sharp} & a
\end{bmatrix}\right).
$$
Therefore the proof is complete.
\end{proof}

%
%

\subsection{Bijection between $\mathcal{L}/\sim$ and
$\mathcal{M}/\sim$}

In this subsection, we shall show that
$\mathcal{F}\circ\mathcal{G} = id_{\mathcal{M}/\sim}$
and
$\mathcal{G}\circ\mathcal{F} = id_{\mathcal{L}/\sim}$.

\begin{lemma}\label{lemma:bimap}
Let $(B,E)$ be an element in $\mathcal{L}$ and $C^*\la B, e_A \ra$
the basic construction for $(B,E)$.
Then for each $x \in C^*\la B, e_A \ra$, there uniquely exists $b
\in B$ such that
$e_Ax = e_Ab$.
\end{lemma}
\begin{proof}
Let $x = \sum_i b_ie_Ac_i$, where $b_i, c_i \in B$.
Then
$$
e_Ax = \sum_i e_Ab_ie_Ac_i
= \sum_i e_AE(b_i)c_i = e_A\sum_i E(b_i)c_i.
$$
And hence $b = \sum_i E(b_i)c_i$.
If
$e_Ab = e_Ab' $, where $b, b' \in B$, then
$$
b = \frac{1}{2}\widetilde{E}(e_A b) =  \frac{1}{2}\widetilde{E}(e_A b') = b',
$$
where $ \widetilde{E}$ is the dual conditional expectation of $C^*\la B, e_A
\ra$
onto $B$.
Thus we obtain the conclusion.
\end{proof}

Let $(B,E)$ be an element in $\mathcal{L}$.
Let $B_{-}$ be a linear subspace of $B$ defined by
$$
 B_{-} = \{ \ b \in B \ | \ E(b) = 0\ \} = \{ \ b \in B \ | \ \beta(b) = -b\
\},
$$
where $\beta$ is an automorphism of $B$ defined in Proposition
\ref{prop:aut}.
By a routine computation we can see that
$B_{-}$ is  an element in $\mathcal{M}$ with
the involution $x^\sharp = x^*$ and the
left and the right $A$-valued inner products defined by
$$
{}_A\la x, y \ra = E(x y^* ) \quad
\la x, y \ra_A = E(x^*y) , \quad \text{for $x, y \in B_{-}$}.
$$
\begin{lemma}\label{lemma:biiso}
With the above notations,
$B_{-} \sim X_B$ i.e., $[B_{-}]= [X_B]$ in $\mathcal{M}/\sim$.
\end{lemma}
\begin{proof}
By Lemma \ref{lemma:bimap}, we can define a map $\varphi$ from
$C^*\la B, e_A \ra$ to $B$ by $e_Ax = e_A\varphi(x)$.
For $e_Ax(1-e_A) \in X_B$, we have
$$
e_Ax(1-e_A)
= e_A\varphi(x)-e_A E(\varphi(x))
= e_A(\varphi(x)-E(\varphi(x))) .
$$
And hence
$$
\varphi(e_Ax(1-e_A)) = \varphi(x)-E(\varphi(x)) \in B_{-}.
$$
It is easy to see that $\varphi|_{X_B}$ is  an $A$-$A$-bimodule isomorphism
of
$X_B$ onto $ B_{-}$.
Furthermore
for $e_Ax(1-e_A), e_Ay(1-e_A)  \in X_B$,
\begin{align*}
{}_A\la e_Ax(1-e_A),e_A y(1-e_A) \ra
&= \psi^{-1}(E((\varphi(x)-E(\varphi(x)))(\varphi(y)-E(\varphi(y)))^*)e_A)
\\
&= E((\varphi(x)-E(\varphi(x)))(\varphi(y)-E(\varphi(y)))^*) \\
&= {}_A\la \varphi(x)-E(\varphi(x)) , \varphi(y)-E(\varphi(y)) \ra.
\end{align*}
Similarly,
$$
\la e_Ax(1-e_A),e_A y(1-e_A) \ra_A
= \la \varphi(x)-E(\varphi(x)) , \varphi(y)-E(\varphi(y)) \ra_A.
$$
And
\begin{align*}
\varphi((e_Ax(1-e_A))^\sharp)
&=\varphi(\widehat{\beta}(e_Ax(1-e_A))^*)
=\varphi(\widehat{\beta}((1-e_A)\varphi(x)^*e_A)) \\
&=\varphi(e_A\varphi(x)^*(1-e_A))
=\varphi(x)^*-E(\varphi(x)^*) \\
&=(\varphi(x)-E(\varphi(x)))^*
= \varphi(e_Ax(1-e_A))^*.
\end{align*}
Hence $X_B \sim B_{-}$ in $\mathcal{M}$.
\end{proof}

\begin{lemma}\label{lemma:const1}
$\mathcal{G}\circ\mathcal{F} = id_{\mathcal{L}/\sim}$.
\end{lemma}
\begin{proof}
For $(B,E) \in \mathcal{L}$, it is easy to see that
$\mathcal{G}([B_{-}]) = [B,E]$.
Since $[X_B] = [B_{-}]$ by the previous lemma,
$\mathcal{G}\circ\mathcal{F}([B,E]) = \mathcal{G}([X_B])
= [B,E]$.
Thus the lemma is proved.
\end{proof}
\begin{lemma}\label{lemma:const2}
$\mathcal{F}\circ\mathcal{G} = id_{\mathcal{M}/\sim}$.
\end{lemma}
\begin{proof}
For $X \in \mathcal{M}$,
$$
(B_X)_{-}
= \{ \ x \in B_X \quad | \quad E_X(x) = 0 \ \}
= \left\{\
\begin{bmatrix}
0 & x \\ \widetilde{x^\sharp} & 0
\end{bmatrix}
\quad | \quad x \in X
 \right\} .
$$
So it is easy to see that $[(B_X)_{-}] = [X]$.
And hence by Lemma \ref{lemma:biiso}
$$
\mathcal{F}\circ\mathcal{G}([X]) = \mathcal{F}([B_X,E_X]) = [(B_X)_{-}] =
[X].
$$
This completes the proof.
\end{proof}

\begin{thm}\label{th:one to one}
There is a $1$-$1$ correspondence between $\mathcal{L}/\sim$
and $\mathcal{M}/\sim$.
\end{thm}
\begin{proof}
This is immediate by Lemmas \ref{lemma:const1} and  \ref{lemma:const2}.
\end{proof}

%
%

\section{Applications}

\subsection{Construction of involutive equivalence bimodules by
$2\Bbb Z$-inner $C^*$-dynamical systems.}

Let $A$ be a unital $C^*$-algebra and $(A, \Bbb Z, \alpha)$ a
$2\Bbb Z$-inner $C^*$-dynamical system which means that
$(A, \Bbb{Z}, \alpha)$ is a $C^* $-dynamical system and that there is
a unitary element $z\in A$ with $\alpha(z)=z$ and 
$\alpha^2 =Ad(z)$. In this case, we can form the restricted crossed
product $A\rtimes_{\alpha/2\Bbb{Z}} \Bbb{Z}$ in the sense of
P. Green\cite{Green:restrict}.
Let $X_{\alpha}$ be the vector space $A$ with the obvious left action of
$A$ on $X_{\alpha}$ and the obvious left $A$-valued inner product, but
we define the right action of $A$ on $X_{\alpha}$ by
$x \cdot a = x\alpha(a)$ for any $x \in X_{\alpha}$ and $a \in A$, and
the right $A$-valued inner product by $\la x, y \ra_A = \alpha^{-1}(x^*y)$ for
any $x, y \in  X_{\alpha}$.

\begin{lemma}
We can define an involution $x \mapsto x^\sharp$ on $X_{\alpha}$ by
$$
x^\sharp = \alpha(x^*)z,
$$
where $z$ is a unitary element of $A$ with
$\alpha(z) = z$ and $\alpha^2 = {\rm Ad}(z)$.
\end{lemma}
\begin{proof}
Since $\alpha(z) = z$ and $\alpha^2 = {\rm Ad}(z)$,
by routine computations, we can see that the map
$x \mapsto x^\sharp$ defined by $x^{\sharp}=\alpha(x^* )z$
is an involution on $X_{\alpha}$.
\end{proof}

\begin{prop}
With the above notations, we suppose that $A$ is simple.
Let $B_{X_{\alpha}}$ be a $C^*$-algebra defined by
$X_{\alpha}$ and
$L$ the linking algebra for $X_{\alpha}$ defined in Section 2. Then the
following conditions are equivalent:
\begin{enumerate}
\item $B_{X_{\alpha}}$ is simple,
\item $A'\cap B_{X_{\alpha}} = \Bbb C \cdot 1$,
\item $B_{X_{\alpha}}'\cap L = \Bbb C \cdot 1$,
\item $\alpha$ is an outer automorphism of $A$.
\end{enumerate}
\end{prop}

\begin{proof}
$(1) \Rightarrow (2)$: By Proposition \ref{prop:aut},
$B_{X_{\alpha}}^{\beta}=A$.
Since $A$ is simple, by Pedersen \cite[Proposition 8.10.12]{Pedersen},
$\beta$ is outer. Hence by Pedersen \cite[Proposition 8.10.13]{Pedersen},
$A' \cap B_{X_{\alpha}} = \Bbb C \cdot 1$.

$(2) \Leftrightarrow (3)$: By Watatani \cite[Proposition
2.7.3]{Watatani:index},
$A' \cap B_{X_{\alpha}}$ is anti-isomorphic to $B_{X_{\alpha}}' \cap
C^*\langle B_{X_{\alpha}}, e_A \rangle$.
This implies the conclusion.

$(2) \Rightarrow (4)$: We suppose that there is a unitary element $w \in A$
such that
$ \alpha = Ad(w)$ Then for any $a \in A$
$$
w \cdot a = w \alpha(a) = aw = a \cdot w.
$$
So it is easy to see that
$$
\begin{bmatrix}
0 & w \\ \widetilde{w^{\sharp}} & 0
\end{bmatrix}
\in A' \cap B_{X_{\alpha}}.
$$
This is a contradiction. Thus $\alpha$ is outer.

$(4) \Rightarrow (1)$: We can identify $L$ with the $C^*$-basic constraction
of
$ A \subset B_{X_{\alpha}}$ by Remark \ref{remark:basic constraction}.
Let $\beta$ be an automorphism of $B_{X_{\alpha}}$ defined
in the same way as in Proposition \ref{prop:aut} and let $\hat{\beta}$ be its dual automorphism.
Then $L^{\hat{\beta}} = B_{X_{\alpha}}$ by Lemma \ref{lem:fix}.
We suppose that $\hat{\beta}$ is inner.
Then there is a unitary element $w =
\begin{bmatrix}
a & x \\ \widetilde{y} & b
\end{bmatrix}
\in L$ such that $\hat{\beta} = Ad(w)$.
Hence for any $c \in A$
$$
\hat{\beta}\left(
\begin{bmatrix}
c & 0 \\ 0 & 0
\end{bmatrix}
\right)
=
\begin{bmatrix} a & x \\ \widetilde{y} & b \end{bmatrix}
\begin{bmatrix} c & 0 \\ 0 & 0 \end{bmatrix}
\begin{bmatrix} a & x \\ \widetilde{y} & b \end{bmatrix} ^*.
$$
Hence we obtain that
$$
\begin{bmatrix} 0 & 0 \\ 0 & c \end{bmatrix} =
\begin{bmatrix} aca^* & ac\cdot y \\
\widetilde{ac^* \cdot y} & \langle c^* \cdot y,
y\rangle_A \end{bmatrix}
$$
for any $c \in A$. Put $c=1$. Then $a= 0$ and $\la y, y \ra_{A}=1$.
Since $w$ is a unitary element, by a routine computation we can see that
$b=0$ and
${}_A \langle y, y \rangle = 1$. This implies that $y$ is
a unitary element in
$A$.
Since $c = \langle c^* \cdot y,y \rangle_A
= \alpha(y^*cy) = \alpha(y)^*\alpha(c)
\alpha(y)$ for any
$c \in A$, $\alpha$ is inner. This is a contradiction. Hence
$\hat{\beta}$ is outer. Since $L$ and $A$ are stably isomorphic by
Brown, Green and Rieffel \cite{BGR:stable},
$L$ is simple. By Pedersen \cite[Theorem 8.10.12]{Pedersen}, $B_{X_{\alpha}}
= L^{\hat{\beta}}$
is simple.
\end{proof}

\begin{lemma}\label{lemma:inner}
Let $(A, \Bbb Z, \alpha)$ be a $2\Bbb Z$-inner dynamical
system with $\alpha(z)=z$ and $\alpha^2 =Ad(z)$, where $z$ is a unitary element in $A$.
Let $B$ is the restricted crossed product $A \rtimes_{\alpha/2\Bbb Z}\Bbb Z$
associated with $(A, \Bbb Z, \alpha)$ and $E$ the canonical conditional
expectation of $B$ onto $A$.
Then $X_B \cong X_{\alpha}$ as involutive $A$-$A$-equivalence bimodules,
where
$X_B$ is an involutive $A$-$A$-equivalence bimodule
induced by $(B, E)$.
\end{lemma}
\begin{proof}
We may assume that $A$ acts on a Hilbert space $H$.
By Olesen and Pedersen \cite[Proposition 3.2]{OP}, we also
assume that $B$ acts on the induced Hilbert space $Ind_{2\Bbb Z}^{\Bbb Z}(H)$.
Let
$$
C=\{\begin{bmatrix} a & x \\
\alpha(xz) & \alpha(a) \end{bmatrix}
\in M_2 (A) \, | \, a,x\in A \} .
$$
Since $A$ acts on $H$, we can $C$ as a $C^{*}$-algebra acting on $H\oplus H$.
We claim that $B\cong C$. Indeed, let $\rho$ be a map
from $K(\Bbb Z , A, z)$ to $C$ defined by for any $f\in K(\Bbb Z , A, z)$
$$
\rho(f)=\begin{bmatrix} f(0) & f(1) \\
\alpha(f(1)z) & \alpha(f(0)) \end{bmatrix} ,
$$
where $K(\Bbb Z ,  A, z)$ is a *-algebra of all
functions $f:\Bbb Z \longrightarrow A$ satisfying that
$f(n-2m)=f(n)z^m $ for any $m,n\in \Bbb Z$ (see Olesen and
Pedersen \cite{OP}). Then by routine computations $\rho$ is
a homomorphism of $K(\Bbb Z , A, z)$ to $C$. Let $U$ be a map
from $Ind_{2\Bbb Z}^{\Bbb Z}(H)$ to $H\oplus H$ defined by
$U\xi=\xi(0)\oplus\xi(1)$ for any $\xi\in K(\Bbb Z , A, z)$.
Then by an easy computaion $U$ is a unitary operator
of $Ind_{2\Bbb Z}^{\Bbb Z}(H)$ onto $H\oplus H$.
Moreover, for any $f\in K(\Bbb Z , A, z)$,
$\rho(f)=UfU^*$. Hence $\rho$ is an isometry
of $K(\Bbb Z , A,z)$ to $C$ and we can extend
$\rho$ to an isomorphism of $B$ onto $C$ since
$K(\Bbb Z , A, z)$ is dense in $B$. Thus $B\cong C$.
Let $F$ be a linear map of $C$ onto $A$ defined by
$F(\begin{bmatrix} a & x \\
\alpha(xz) & \alpha(a) \end{bmatrix})
=\begin{bmatrix} a & 0 \\
0 & \alpha(a) \end{bmatrix}$ for any
$\begin{bmatrix} a & x \\
\alpha(xz) & \alpha(a) \end{bmatrix}\in C$,
where we identify $A$ with a $C^*$-algebra
$\{\begin{bmatrix} a & 0 \\
0  & \alpha(a) \end{bmatrix} \, | \, a\in A \}$.
Then by an easy computation $(B, E)\sim (C, F)$ in $\mathcal{L}$.
Let $(B_{X_{\alpha}}, E_{X_{\alpha}})$ be
an element in $\mathcal{L}$ induced by the involutive
$A$-$A$-equivalence bimodule $X_{\alpha}$. Let $\Phi$ be
a map from $C$ to $B_{X_{\alpha}}$ defined by
$$
\Phi(\begin{bmatrix} a & x \\
\alpha(xz) & \alpha(a) \end{bmatrix})
=\begin{bmatrix} a & x \\
\widetilde{x^{\sharp}} & a \end{bmatrix}
$$
for any $\begin{bmatrix} a & x \\
\alpha(xz) & \alpha(a) \end{bmatrix}\in C$.
Then by routine computations $\Phi$ is an
isomorphism of $C$ onto $B_{X_{\alpha}}$ with
$F=E_{X_{\alpha}}\circ\Phi$.
Thus $(B, E)\sim(B_{X_{\alpha}}, E_{X_{\alpha}})$.
By Theorem \ref{th:one to one}, $X_B \sim X_{\alpha}$ in $\mathcal{M}$.
\end{proof}

Let $B$ be a unital $C^*$-algebra and $A$ a $C^*$-subalgebra of $B$ with a
common unit.
Let $E$ be a conditional expectation of $B$ onto $A$ with $ \ind E = 2$.
For any $n \in \Bbb N$ let $M_n$ be the $n \times n$-matrix algebra over
$\Bbb C$ and $M_n(A)$ the $n \times n$-matrix algebra over $A$.
Let $\{(x_i, x_i^*)\}_{i=1}^n$ be a quasi-basis for $E$.
We define $q = [q_{ij}] \in M_n(A)$ by $q_{ij} = E(x_i^*x_j)$.
Then by Watatani \cite{Watatani:index}, $q$ is a projection
and $C^*\la B, e_A\ra \simeq qM_n(A)q$.
Let $\pi$ be an isomorphism of $C^*\la B, e_A\ra$ onto
$qM_n(A)q$ defined by
$$
\pi(ae_Ab) = [E(x_i^*a)E(bx_j)] \in M_n(A)
$$
for any $a,b \in B$.
Especially for any $b \in B$,
$\pi(b) =  [E(x_i^*bx_j)]$
since $\sum_{i=1}^nx_ie_Ax_i^*=1$.

\begin{prop}\label{prop:inner sys}
With the above notations,
the following conditions are equivalent:
\begin{enumerate}
\item $e_A$ and $1-e_A$ are equivalent in $C^*\la B, e_A\ra$,
\item there exists a unitary element $u \in B$ such that
$\{(1,1), (u, u^*)\}$ is a quasi-basis for $E$,
\item there exists a $2\Bbb Z$-inner $C^*$-dynamical system $(A, \Bbb Z,
\alpha)$
such that $X_{\alpha} \sim X_B$.
\end{enumerate}
\end{prop}
\begin{proof}
$(1) \Rightarrow (2)$: We suppose that there is a partial
isometry $v \in C^*\la B, e_A \ra$ such that
$v^*v = e_A $, $\ vv^* = 1-e_A$.
Then
$ve_Av^* = 1-e_A $.
By Lemma \ref{lemma:bimap}, there exists an element $u$ in $B$ such that
$ve_A = ue_A$ and hence
$ue_Au^* = 1-e_A$.
Let $\widetilde{E}$ be the dual conditional expectation for $E$.
Then
$$
uu^* = 2 \widetilde{E}(ue_Au^*) = 2\widetilde{E}(1-e_A) =1.
$$
Therefore $u$ is a co-isometry element in $B$.
Since $e_A u^*u e_A = e_Av^*ve_A = e_A$, we have
$E(u^*u) = 1$ and $E(1- u^*u) = 0$.
And hence $u^*u = 1$  i.e., $u$ is a unitary element in $B$.
For any $x \in B$
$$
xe_A
= (e_A + ue_Au^*)xe_A
= E(x)e_A + uE(u^*x)e_A
= (E(x) + uE(u^*x))e_A.
$$
Thus $x = E(x) + uE(u^*x)$ by Lemma \ref{lemma:bimap}.
Similarly, $x= E(x) + E(xu)u^*$.
This implies that $\{(1,1), (u, u^*)\}$ is a quasi-basis for $E$.

$(2) \Rightarrow (1)$:
We suppose that $\{(1,1), (u, u^*)\}$ is a quasi-basis for $E$ and that $u$ is a
unitary element in $B$.
Then
$$
u = E(u) + uE(u^*u) = E(u) + u.
$$
This implies that $E(u)=0$. Hence
$$
q =
\begin{bmatrix}
E(1\cdot 1) & E(u) \\
E(u^*) & E(u^*u)
\end{bmatrix}
=
\begin{bmatrix}
1 & 0 \\
0 & 1
\end{bmatrix}.
$$
Therefore $C^*\la B, e_A \ra \simeq M_2(A)$.
Furthermore
$$
\pi(e_A) =
\begin{bmatrix}
1 & 0 \\
0 & 0
\end{bmatrix}, \quad
\pi(1-e_A) =
\begin{bmatrix}
0 & 0 \\
0 & 1
\end{bmatrix}.
$$
And hence $e_A \sim (1-e_A)$ in $C^*\la B, e_A \ra$.

$(2) \Rightarrow (3)$:
We suppose that $\{(1,1), (u, u^*)\}$ is a quasi-basis for $E$ and that $u$ is a
unitary element in $B$.
Then in the same way as above $E(u)=0$.
For any $a \in A$
$$
uau^* = E(uau^*) + E(uau^*u)u^* = E(uau^*) +E(u)au^* = E(uau^*).
$$
Therefore $uAu^* = A$. Let $\alpha$ be an automorphism of $A$ defined by
$\alpha(a) =  uau^*$ for any $a\in A$.
Since $u^2 = E(u^2) + uE(u^*u^2) = E(u^2)$,  $ u^2$ is an element in $A$.
Therefore $(A, \Bbb Z, \alpha)$ is a $2\Bbb Z$-inner $C^*$-dynamical system.
It is easy to see that
$$
X_{\alpha} \sim Au = B_{-} = \{ \ b \in B \  | \ E(b) = 0\ \}.
$$
By Lemma \ref{lemma:biiso}, $X_{\alpha} \sim X_B$.

$(3) \Rightarrow (2)$:
We suppose that there exists a $2\Bbb Z$-inner $C^*$-dynamical system $(A,
\Bbb Z, \alpha)$
such that $X_{\alpha} \sim X_B$.
By the previous lemma, we may suppose that
$B =  A \rtimes_{\alpha/2\Bbb Z}\Bbb Z$.
Then there exists a unitary element $u \in B$ such that
${\rm Ad}(u) = \alpha$, $u^2 \in A$ and $E(u) = 0$.
By a routine computation we can see that $\{(1,1), (u, u^*)\}$ is
a quasi-basis for $E$.
\end{proof}

\begin{cor}\label{cor:irrational}
Let $\theta$ be an irrational number in $(0,1)$ and $A_{\theta}$ the
corresponding irrational rotation
$C^*$-algebra. Let $B$ be a unital $C^*$-algebra including $A_{\theta}$ as a
$C^*$-subalgebra of $B$ with a common unit.
We suppose that there is a conditional expectation $E$ of $B$ onto
$A_{\theta}$ with
$\ind E =2 $.
Then there is a $2\Bbb Z$-inner $C^*$-dynamical system $(A_{\theta}, \Bbb
Z, \alpha)$
such that $(B,E) \sim (A_{\theta} \rtimes_{\alpha/2\Bbb Z}\Bbb Z, F)$,
where $F$ is the canonical conditional expectation of  $A_{\theta}
\rtimes_{\alpha/2\Bbb Z}\Bbb Z$ onto
$A_{\theta}$.
\end{cor}
\begin{proof}
Let $e$ be the Jones projection induced by $E$.
We can identify the basic construction $C^*\la B, e \ra$ with
$q M_n(A_{\theta})q$ in the same way as in the
previous argument. Hence $C^*\la B, e \ra$ has
the unique normalized
trace $\tau$ and $\tau(e) = \tau(1-e) = \frac{1}{2}$.
So it is easy to see that $ e \sim 1-e $ in $C^*\la B, e \ra$
since $A_{\theta}$ has cancellation.
Therefore we obtain the conclusion by the previous proposition.
\end{proof}

%
%
%

\subsection{Examples}
In this subsection,
let $A_{\theta}$ be as in Corollary\ref{cor:irrational} and
let $u$, $v$ be two unitary generators
satisfying the
commutation relation:
$$
uv = e^{2\pi i \theta}vu.
$$

\begin{exam}
Let $A_{2\theta}$ be the $C^*$-subalgebra of $A_{\theta}$ generated by $
u^2$ and $v$.
Then we can denote $A_{\theta} = \{ x+ yu \ | \  x, y \in A_{2\theta} \}$.
Let $E$ be a map of $A_{\theta}$ onto $A_{2\theta}$ defined by $E(x+yu) =
x$.
It is easy to see that $E$ is a conditional expectation of $A_{\theta}$
onto $A_{2\theta}$ with
$\ind E = 2$ and a quasi-basis $\{(1,1),(u, u^*)\}$. Hence by Corollary
\ref{cor:irrational},
$A_{\theta}$ can be represented as the restricted crossed product $A_{2\theta}
\rtimes_{\alpha/2\Bbb Z} \Bbb Z$,
where $\alpha$ is an automorphism on $A_{2\theta}$ defined by $\alpha =
Ad(u)$.

Suppose that $A_{\theta}$ can be represented as a crossed product
$A_{2\theta}\rtimes_{\beta}\mathbb{Z}_2 $
for some $\mathbb{Z}_2 $-action $\beta$ on $A_{2\theta}$.
Then there exists a
self-adjoint unitary element $w$ in $A_{\theta}$ satisfying that
$\beta = Ad(w)$ and $A_{\theta} = \{x + yw \ | \ x, y \in A_{2\theta} \}$.
Let $\tau$ be the unique tracial state on $A_{\theta}$. By the uniqueness of
$\tau$, we can see that
$\tau(x+yw) = \tau (x)$.  Let $e$ be a projection in $A_{\theta}$ defined by
$ e= \frac{1}{2}(1 + w)$.
Then $\tau(e) = \frac{1}{2}$. This contradicts that $\tau(A_{\theta}) =
\left(\Bbb Z \cap \theta \Bbb Z \right) \cap (0, 1)$. Therefore $A_{\theta}$
can not be represented as
a crossed product $A_{2\theta}\rtimes_{\beta}\mathbb{Z}_2 $
for any $\mathbb{Z}_2 $-action $\beta$ on $A_{2\theta}$.
\end{exam}

\begin{exam}
Let $\sigma$ be the involutive automorphism of $A_{\theta}$ determined by
$\sigma(u) = u^*$ and
$\sigma(v) = v^*$. Let $C_{\theta}$ denote the fixed point algebra
$A_{\theta}^{\sigma} =
\{ x \in A_{\theta} \ | \ \sigma(x) = x \}$ and $B_{\theta}$ the crossed
product
$A_{\theta} \rtimes_{\sigma}\mathbb{Z}_2 $. Then $B_{\theta}$ is the basic
construction of
$C_{\theta} \subset A_{\theta}$.
By Kumjian\cite{Kum1}, $K_0$-group of $B_{\theta}$,
$K_0 (B_{\theta})$ is isomorphic to
${\Bbb Z}^6$.
By routine computations, we can see $[e] \not= [1-e]$ in $K_0(B_{\theta})$,
where
$e$ is the Jones projection for the inclusion $ C_{\theta}\subset
A_{\theta}$.
Hence $e \not\sim 1-e$ in $B_{\theta}$. Therefore the inclusion $C_{\theta}
\subset A_{\theta}$
can not be represented as the restricted crossed product
$C_{\theta} \subset C_{\theta} \rtimes_{\alpha/2\Bbb Z} \Bbb Z$ for any
automorphism
$\alpha$ on  $C_{\theta}$ by Proposition \ref{prop:inner sys}.

\end{exam}

\bf
Acknowledgement.
\rm
The authors wish to thank the referee for some
valuable suggestions for improvement of the manuscript,
especially for improvement of Introduction.


\begin{thebibliography}{99}


\bibitem{BEEK1}
O.\ Bratteli, G.\ A,\ Elliott, D.\ E.\ Evans, A.\ Kishimoto,
{\it Non-commutative spheres I},
Int.\ J.\ Math.\ ,
{\bf 2}
(1991),
p.\ 139--166.

\bibitem{BGR:stable}
L.\ G.\ Brown, P.\ Green and M.\ A.\ Rieffel,
{\it Stable isomorphism and strong Morita equivalence of
$C^*$-algebra},
Pacific J.\ Math.\ ,
{\bf 71}
(1977),
p.\ 349--368.

\bibitem{ER}
G.\ A.\ Elliott and M.\ R\o rdam,
{\it The automorphism group of the irrational rotation algebra},
Comm.\ Math.\ Phys.\ ,
{\bf 155}(1993),
p.\ 3--26.



\bibitem{Green:restrict}
P.\ Green,
{\it The local structure of twisted covariance algebras},
Acta Math.\ ,
{\bf 140}(1978),
p.\ 191--250.

\bibitem{Izumi:subal}M.\ Izumi,
{\it Inclusions of simple $C^*$-algebras},
J.\ reine angew.\ Math. \ ,
{\bf 547}(2002),
p.\ 97--138.

\bibitem{Jones:factor}
V.\ Jones,
{\it Index for subfactors},
Invent.\ Math.,
{\bf 72}(1983), p.\ 1--25.

\bibitem{KW}
T.\ Kajiwara and Y.\ Watatani,
{\it Jones index theory by Hilbert $C^*$-bimodules and $K$-Theory},
Trans.\ Amer.\ Math.\ Soc.\ ,
{\bf 352}(2000), p.\ 3429--3472.

\bibitem{Kum1}
A.\ Kumjian,
{\it On the K-theory of the symmetrized non-commutative torus},
C.\ R.\ Math.\ Rep.\ Acad.\ Sci.\ Canada,
{\bf 12}(1990),
p.\ 87--89.

\bibitem{OP}
D.\ Olesen and G.\ K.\ Pedersen,
{\it Partially inner $C^*$-dynamical systems},
J.\ Funct.\ Anal.\ ,
{\bf 66}(1986),p.\ 262--281.

\bibitem{Pedersen}
G.\ K.\ Pedersen,
{$C^*$-algebras and their automorphism groups},
Academic Press,
1979.

\bibitem{Reiffel:irr}
M.\ A.\ Rieffel,
{\it $C^*$-algebra associated with irrational rotations},
Pacific J.\ Math.\,
{\bf 93}(1981),
p.\ 415--429.

\bibitem{Watatani:index} Y.\ Watatani,
{\it Index for $C^*$-subalgebras},
Mem.\ Amer.\ Math.\ Soc.\ ,
{\bf 424},
Amer.\ Math.\ Soc., \ Providence, R.\ I.,
(1990).
\end{thebibliography}
\end{document}